\documentclass[12pt]{amsart}
\usepackage{amssymb,epsfig,times}

\newtheorem{Thm}{Theorem}
\newtheorem{Cor}[Thm]{Corollary}
\newtheorem{Lemma}[Thm]{Lemma}
\newtheorem{Prop}[Thm]{Proposition}

\theoremstyle{definition}
\newtheorem{Defn}{Definition}
\newtheorem{Notation}[Defn]{Notation}
\newtheorem{Remark}[Thm]{Remark}
\newtheorem{Ex}[Thm]{Example}

\newcommand{\norm}[1]{\left\Vert#1\right\Vert}
\newcommand{\abs}[1]{\left\vert#1\right\vert}
\newcommand{\chf}[1]{\mathbf{1}_{#1}}

\newcommand{\set}[1]{\left\{#1\right\}}

\DeclareMathOperator{\Part}{\mathcal{P}}
\newcommand{\St}[1]{\mathrm{St}_{#1}}
\DeclareMathOperator{\NC}{\mathit{NC}}

\newcommand{\mf}[1]{\mathbb{#1}}
\newcommand{\mc}[1]{\mathcal{#1}}
\newcommand{\mb}[1]{\mathbf{#1}}

\DeclareMathOperator{\Sym}{Sym}
\newcommand{\id}{\mathrm{Id}}
\newcommand{\ip}[2]{\left \langle #1, #2 \right \rangle}
\newcommand{\state}[1]{\varphi \left[ #1 \right]}

\DeclareMathOperator{\Falg}{\mathcal{F}_{\mathrm{alg}}}
\newcommand{\rc}[1]{\mathrm{rc} \left( #1 \right)}
\newcommand{\ls}[1]{\mathrm{span} \left( #1 \right)}
\newcommand{\Exp}[1]{\left\langle #1 \right\rangle}

\DeclareMathOperator{\E}{\mathrm{E}}
\renewcommand{\phi}{\varphi}
\newcommand{\Aalg}{\mathcal{A}^{\mathrm{alg}}}

\allowdisplaybreaks[1]

\title[]{$q$-L\'{e}vy processes}
\author[M.~Anshelevich]{Michael Anshelevich}
\thanks{This work was supported in part by an NSF postdoctoral fellowship}
\address{Department of Mathematics, University of California, Riverside, CA 92521-0135}
\email{manshel@math.ucr.edu}
\date{\today}

\begin{document}

\begin{abstract}
We continue the investigation of the L\'{e}vy processes on a $q$-deformed full Fock space started in \cite{AnsQCum}. First, we show that the vacuum vector is cyclic and separating for the algebra generated by such a process. Next, we describe a chaotic representation property for it in terms of multiple integrals with respect to diagonal measures, in the style of Nualart and Schoutens. We define stochastic integration with respect to these processes, and calculate their combinatorial stochastic measures. Finally, we show that they generate infinite von Neumann algebras.
\end{abstract}

\maketitle

\section{Introduction}
The original motivation for this paper was to understand the results about a chaos decomposition for general L\'{e}vy processes obtained by Nualart and Schoutens \cite{SchChaotic}. A relation between multiple Wiener integrals for the Brownian motion and the symmetric Fock space, and its role in the chaos decomposition, are well known. For more general L\'{e}vy processes, a representation on a symmetric Fock space was introduced by Parthasarathy, see \cite{Parthasarathy}. We will show that in this representation, the chaos decomposition has a very natural interpretation, and its proof becomes more immediate. A recent preprint \cite{Vershik-Levy}, which came to our attention when the paper was nearing completion, handles a number of related questions, although by quite different methods. A number of preceding papers which deal with related topics are cited in the references.

Our second goal was the extension of these results to certain non-commutative stochastic processes. Here the starting point is the $q$-Fock space construction introduced by Bo\.{z}ejko and Speicher \cite{BozSpeBM1,BKSQGauss}. They also introduced the $q$-Brownian motion, which reduces to the usual Brownian motion for the bosonic case $q=1$, to the fermionic Brownian motion for $q=-1$, and to the analog of the Brownian motion in free probability for $q=0$. In a previous paper \cite{AnsQCum}, following the ideas of Parthasarathy and Sch\"{u}rmann \cite{SchurCondPos}, we introduced L\'{e}vy processes on the $q$-deformed Fock spaces. Since we are able to express the chaos decomposition property in the Fock space language, the same method gives the corresponding result for the $q$-deformed processes. Moreover, this decomposition can then be used to obtain a number of consequences.

Combinatorial stochastic measures were introduced for the usual stochastic processes by Rota and Wallstrom \cite{Rota}. In this paper, we obtain an explicit formula for such measures for the $q$-L\'{e}vy processes, unifying a number of previous results. $k$'th power of a process can be expressed as a sum over stochastic measures. On the level of operators, this turns out to be precisely the decomposition of the product of operators in terms of their Wick products.

Since the $q$-L\'{e}vy processes consist of non-commuting operators, they generate non-commutative von Neumann algebras, and any information about these algebras is of interest. It was known that for the $q$-Brownian motion, or for the free probability case $q=0$, these algebras are II$_1$-factors. We show that in contrast, in all other cases these algebras have no normal tracial states, and thus are infinite algebras. Further information about these algebras awaits discovery.

There are three standard classes of infinitely divisible distributions; see \cite{Streater} for one interpretation of these classes. The intermediate one, with Fourier transforms
\[
\mc{F}_\mu(\theta) = \exp \left[ \int_{\mf{R}} (e^{i \theta x} - 1 - i \theta x) \frac{1}{x^2} \,d\nu(x) \right]
\]
for $\nu$ a finite measure, was considered by Kolmogorov. It consists of the distributions with mean $0$ and finite variance. This is the class whose $q$-analogs are treated in this paper. For simplicity we treat only one-dimensional such processes, although there is no real difficulty in extending the results to $k$ dimensions.

The narrowest class, with Fourier transforms
\[
\mc{F}_\mu(\theta) = \exp \left[ \int_{\mf{R}} (e^{i \theta x} - 1) \,d\nu(x) \right],
\]
for $\nu$ a finite measure, consists of the compound Poisson distributions, considered by de Finetti. In this case there is a more natural construction in which the starting object is an algebra with a fixed state. The analogs of some results of the paper for this class are treated in the appendix. Finally, all infinitely divisible distributions have Fourier transforms of the form
\[
\mc{F}_\mu(\theta) 
= \exp \left[ i \gamma \theta + \int_{\mf{R}} \left( e^{i \theta x} - 1 - \frac{i \theta x}{1 + x^2} \right) \frac{1 + x^2}{x^2} \,d\nu(x) \right]
\]
for $\nu$ a finite measure, by the L\'{e}vy-Khinchine theorem. Representations of all of the corresponding L\'{e}vy processes on the symmetric Fock space are described in Section 21 of~\cite{Parthasarathy}. There is no difficulty with extending the definitions to the general $q$ case, but since the usual treatment involves measures which are no longer finite, we do not treat this case in this paper.

\noindent\textbf{Acknowledgments:} I thank Ed Effros, Marius Junge, and Murad Taqqu for a number of useful (and enjoyable) conversations.

\section{Operators on the $q$-Fock space}
This paper is a sequel to \cite{AnsQCum}; see that paper for all definitions and references not explicitly provided here.

\subsection{$q$-deformed full Fock space}
\label{Sec:q-Fock}
Let $\mc{H}_0$ be a real Hilbert space and $\mc{H}$ its complexification. Let $\mc{H}^{\otimes n}$, for $n \geq 0$, be its ``$n$-particle space.'' Let its algebraic Fock space be the vector space
\[
\Falg(\mc{H}) 
= \bigoplus_{n=0}^\infty \mc{H}^{\otimes n}
= (\mf{C} \Omega) \oplus \mc{H} \oplus \mc{H}^{\otimes 2} \oplus \mc{H}^{\otimes 3} \oplus \ldots
\]
Here $\Omega$ is the generator of the $0$'th component, traditionally called the vacuum vector. Define an inner product on $\Falg(\mc{H})$ by
\[
\ip{\xi_1 \otimes \xi_2 \otimes \ldots \otimes \xi_n}{\eta_1 \otimes \eta_2 \otimes \ldots \otimes \eta_k}_0
= \delta_{n k} \ip{\xi_1}{\eta_1} \ip{\xi_2}{\eta_2} \ldots \ip{\xi_n}{\eta_n},
\]
where $\ip{\cdot}{\cdot}$ is the inner product on $\mc{H}$.

Fix $q \in (-1,1)$. Define an operator $P_n$ on $\mc{H}^{\otimes n}$ by
\[
P_n(\xi_1 \otimes \xi_2 \otimes \ldots \otimes \xi_n)
= \sum_{\sigma \in \Sym(n)} q^{i(\sigma)} \xi_{\sigma(1)} \otimes \xi_{\sigma(2)} \otimes \ldots \otimes  \xi_{\sigma(n)},
\]
where $i(\sigma)$ is the number of inversions of the permutation $\sigma$. Then all $P_n$, and so the corresponding operator $P$ on the whole of $\Falg(\mc{H})$, are strictly positive \cite{BozSpeBM1}. We define a new $q$-inner product on $\Falg(\mc{H})$ by
\[
\ip{\vec{\xi}}{\vec{\eta}}_q
= \ip{\vec{\xi}}{P \vec{\eta}}_0.
\]
This is a positive definite inner product. Completing the algebraic Fock space with respect to the corresponding norm, we get the $q$-Fock space $\mc{F}_q(\mc{H})$.

For $q=\pm 1$, the inner product is only positive semi-definite. For $q=1$, taking the quotient by its kernel gives the symmetric Fock space; for $q=-1$, it gives the anti-symmetric Fock space. For $q=0$, we get the full Fock space.

For $\zeta \in \mc{H}_0$, define creation and annihilation operators on $\mc{F}_q(\mc{H})$ by
\begin{align*}
& a^\ast(\zeta) \Omega = \zeta, \\
& a^\ast(\zeta) (\eta_1 \otimes \ldots \otimes \eta_n) = \zeta \otimes \eta_1 \otimes \ldots \otimes \eta_n, \\
& a(\zeta) \Omega = 0, \\
& a(\zeta) \eta = \ip{\zeta}{\eta} \Omega, \\
& a(\zeta) (\eta_1 \otimes \ldots \otimes \eta_n) = \sum_{k=1}^n q^{k-1} \ip{\zeta}{\eta_k} \eta_1 \otimes \ldots \otimes \check{\eta}_k \otimes \ldots \otimes \eta_n.
\end{align*}
Here $\check{\eta}$ means ``omit $\eta$''. These operators can be extended to bounded operators on the whole $\mc{F}_q(\mc{H})$, so that $a(\zeta)$ and $a^\ast(\zeta)$ are adjoints of each other. Moreover,
\[
a(\zeta) a^\ast(\eta) - q a^\ast(\eta) a(\zeta) = \ip{\zeta}{\eta} \id.
\]
Let $T$ be a bounded operator on $\mc{H}_0$. We will also denote by $T$ its complexification, which is a bounded operator on $\mc{H}$. Define the gauge (or preservation, or differential second quantization) operator $p(T)$ by 
\begin{align*}
&p(T) \Omega = 0, \\
&p(T) (\eta_1 \otimes \ldots \otimes \eta_n)
= \sum_{k=1}^n q^{k-1} (T \eta_k) \otimes \eta_1 \otimes \ldots \otimes \check{\eta}_k \otimes \ldots \otimes \eta_n.
\end{align*}
By Proposition 2.2 of~\cite{AnsQCum}, if $T$ is a self-adjoint operator, the operator $p(T)$ is essentially self-adjoint with dense domain $\Falg(\mc{H})$.

\begin{Lemma}
The operator $p(T)$ on the $q$-Fock space is bounded.
\end{Lemma}

\begin{proof}
Denote by $\norm{p(T)}_{q \mapsto q}$ the norm of $p(T)$ as an operator from $\mc{F}_q(\mc{H})$ to itself. First we show that $\norm{p(T)}_{0 \rightarrow 0} < \infty$. Indeed, $p(T)$ is a composition of two operators,
\[
\pi(\eta_1 \otimes \ldots \otimes \eta_n)
= \sum_{k=1}^n q^{k-1} \eta_k \otimes \eta_1 \otimes \ldots \otimes \check{\eta}_k \otimes \ldots \otimes \eta_n
\]
followed by
\[
p_0(T) (\eta_1 \otimes \ldots \otimes \eta_n)
= (T \eta_1) \otimes \eta_2 \otimes \ldots \otimes \eta_n.
\]
On $\mc{F}_0(\mc{H})$,
\[
\norm{\pi} \leq \max_n \sum_{k=1}^n q^{k-1} \leq \max \Bigl(1, \frac{1}{1-q} \Bigr)
\]
and $\norm{p_0(T)} \leq \norm{T}$ from general tensor product considerations.

Now we show that $\norm{p(T)}_{q \mapsto q} < \infty$. For the remainder of the proof, we work with the $0$-inner product, the $q$-inner product being given by the positive definite density operator $P$. It follows from the proof of Proposition 2.2 of~\cite{AnsQCum} that the adjoint of $p(T)$ with respect to the $q$-inner product is $p(T^\ast)$. This means that $P p(T^\ast) = p(T)^\ast P$, where $^\ast$ denotes the adjoint with respect to the $0$-inner product. In particular,
\[
P p(T^\ast) p(T) = p(T)^\ast P p(T) \geq 0
\]
is a positive operator. So from the inequality
\[
P p(T^\ast) p(T)^2 p(T^\ast) P \leq \norm{p(T^\ast) p(T)^2 p(T^\ast)}_{0 \rightarrow 0} P^2
\]
it follows that 
\[
P p(T^\ast) p(T) \leq \sqrt{\norm{p(T^\ast) p(T)^2 p(T^\ast)}_{0 \rightarrow 0}} P \leq \norm{p(T^\ast)}_{0 \rightarrow 0} \norm{p(T)}_{0 \rightarrow 0} P
\]
So
\begin{align*}
\ip{p(T) \vec{\xi}}{p(T) \vec{\xi}}_q 
&= \ip{\vec{\xi}}{p(T^\ast) p(T) \vec{\xi}}_q
= \ip{\vec{\xi}}{P p(T^\ast) p(T) \vec{\xi}}_0 \\
&\leq \norm{p(T^\ast)}_{0 \rightarrow 0} \norm{p(T)}_{0 \rightarrow 0} \ip{\vec{\xi}}{P \vec{\xi}}_0 \\
&= \norm{p(T^\ast)}_{0 \rightarrow 0} \norm{p(T)}_{0 \rightarrow 0} \ip{\vec{\xi}}{\vec{\xi}}_q.
\end{align*}
Since $\norm{T} = \norm{T^\ast}$, we conclude that
\[
\norm{p(T)}_{q \rightarrow q} \leq \sqrt{\norm{p(T^\ast)}_{0 \rightarrow 0} \norm{p(T)}_{0 \rightarrow 0}} \leq \max \Bigl(1, \frac{1}{1-q} \Bigr) \norm{T}.
\]
\end{proof}

We conclude that for any pair $(\zeta, T)$, $\zeta \in \mc{H}_0$, the operator
\[
p(\zeta, T) = a(\zeta) + a^\ast(\zeta) + p(T)
\]
is a bounded operator on the $q$-Fock space of $\mc{H}$, self-adjoint if $T$ is. Note that $p(\cdot)$ is linear as a function of the pair $(\zeta, T)$.

\subsection{Construction from a generator} 
The following  construction is inspired by \cite{SchurCondPos,GloSchurSpe}. Let $\mc{S}$ be an index set. Denote by $\mc{S}^\infty$ the set of multi-indices (finite sequences) of elements of $\mc{S}$. Denote by $\mf{R}_0 \langle \mb{x}, \mc{S} \rangle$ the algebra of all real polynomials without a constant term, in a collection of non-commuting indeterminates $\set{x_i: i \in \mc{S}}$. It is an algebra without a unit, and even a star-algebra with the obvious involution. Let $\psi$ be a positive linear functional on $\mf{R}_0 \langle \mb{x}, \mc{S} \rangle$; equivalently, it is a conditionally positive linear functional on the full algebra $\mf{R} \langle \mb{x}, \mc{S} \rangle$. Assume that it is also both left- and right-bounded, in the sense that for all $f, g \in \mf{R}_0 \langle \mb{x}, \mc{S} \rangle$ there exist constants $M_f, N_g$ such that 
\[
\psi[f^\ast g^\ast g f] \leq M_f \psi[g^\ast g]; \qquad \psi[f^\ast g^\ast g f] \leq N_g \psi[f^\ast f].
\]
$\psi$ induces a positive semi-definite inner product on the space $\mf{R}_0 \langle \mb{x}, \mc{S} \rangle$ in the usual way, $\ip{f}{g}_{\psi} = \psi[f^\ast g]$, as well as a semi-norm $\norm{\cdot}_\psi$. Taking a quotient by the subspace of seminorm-zero vectors and completing with respect to the induced norm, we obtain a real Hilbert space $\mc{H}_0$ with the induced inner product. Denote by $\rho$ the canonical mapping $\mf{R}_0 \langle \mb{x}, \mc{S} \rangle \rightarrow \mc{H}_0$, let $\mc{D}_0$ be its image, and for $f, g \in \mf{R}_0 \langle \mb{x}, \mc{S} \rangle$ define the operator $\Gamma(f): \mc{D}_0 \rightarrow \mc{D}_0$ by $\Gamma(f) \rho(g) = \rho(f g)$. Put, for $i \in \mc{S}$, $\xi_i = \rho(x_i)$, $T_i = \Gamma(x_i)$. More generally, for $f \in \mf{R}_0 \langle \mb{x}, \mc{S} \rangle$, denote $\xi_f = \rho(f)$, $T_f = \Gamma(f)$. Then each $T_{f}$ is essentially self-adjoint, with dense domain $\mc{D}_0$ consisting of analytic vectors. Since $\psi$ is left-bounded, for each $f \in \mf{R}_0 \langle \mb{x}, \mc{S} \rangle$, $\norm{T_f} \leq N_f$, so each $T_f$ is bounded. Any family $\set{\xi_i, T_i: i \in \mc{S}}$ satisfying a certain compatibility condition arises in this way, see Proposition 4.3 of~\cite{AnsQCum} and Theorem 3 of~\cite{GloSchurSpe}.

Let $\mc{H}$ be the complexification of $\mc{H}_0$ and $\mc{D}$ the complexification of $\mc{D}_0$. Denote $X(i) = p(\xi_i, T_i)$, and more generally $X(f) = p(\xi_f, T_f)$ for $f \in \mf{R}_0 \langle \mb{x}, \mc{S} \rangle$. These are bounded operators on $\mc{F}_q(\mc{H})$. Note that the mapping $f \mapsto X(f)$ is $\mf{R}$-linear. If $\norm{f}_\psi = 0$, by definition $\xi_f = 0$. Moreover, since $\psi$ is right-bounded, $\Gamma(f) = 0$ as well. As a result, $X(f)$ depends only on the equivalence class of $f$ in $\mc{H}$, and so can be defined for $f \in \mc{D}_0$.

Denote by $\Aalg_X$ the complex algebra (no closure) generated by $\set{X(i): i \in \mc{S}}$, and by $\Aalg_{X, \Delta}$ the complex algebra generated by $\set{X(f): f \in \mf{R}_0 \langle \mb{x}, \mc{S} \rangle}$.

\begin{Lemma}
\label{Lemma:Cyclic}
The vacuum vector $\Omega$ is cyclic for $\Aalg_{X, \Delta}$.
\end{Lemma}

\begin{proof}
By definition,
\[
X(f_1) X(f_2) \ldots X(f_n) (\Omega)
= \xi_{f_1} \otimes \xi_{f_2} \otimes \ldots \otimes \xi_{f_n} + \vec{\eta}, 
\]
where $\vec{\eta} \in \bigoplus_{i=0}^{n-1} \mc{H}^{\otimes i}$. It follows by induction that $\Aalg_{X, \Delta} \Omega \supset \Falg(\mc{D})$, which is dense in $\mc{F}_q(\mc{H})$. 
\end{proof}

\begin{Remark}
The result of the preceding lemma will be made more precise through the use of multiple stochastic integrals and Kailath-Segall polynomials.
\end{Remark}

\subsection{The Wick map}
Define the maps $W_0: \mf{R}_0 \langle \mb{x}, \mc{S} \rangle^n \rightarrow \Aalg_{X, \Delta}$, $n = 1, 2, \ldots$, inductively as follows: $W_0(f) = X_{f}$ and
\begin{multline*}
W_0(f_0, f_1, \ldots, f_n) 
= X(f_0) W(f_1, \ldots, f_n) - \sum_{i=1}^n q^{i-1} \ip{\xi_{f_0}}{\xi_{f_i}} W(f_1, \ldots, \check{f}_i, \ldots, f_n) \\
- \sum_{i=1}^n q^{i-1} W(f_0 f_i, f_1, \ldots, \check{f}_i, \ldots, f_n).
\end{multline*}
By the discussion in the preceding subsection, $W_0$ is a multi-linear map which depends only on the equivalence classes of $f_1, f_2, \ldots, f_n$ in $\mc{H}$. So we can project it to a map on $\mc{D}_0^n$, and extend it to a $\mf{C}$-linear map on $\mc{D}^{\otimes n}$. As a result, we can define a linear \emph{Wick map}
\[
W: \Falg(\mc{D}) \rightarrow \Aalg_{X, \Delta}
\]
as such an extension
\[
W(\xi_{f_1} \otimes \ldots \otimes \xi_{f_n})
= W_0(f_1, \ldots, f_n),
\]
with the extra condition $W(\Omega) = \id$. Clearly,
\begin{equation}
\begin{split}
\label{Wick-recursion}
W(\xi_{f_0} \otimes \xi_{f_1} \otimes \ldots \otimes \xi_{f_n}) 
&= X(f_0) W(\xi_{f_1} \otimes \ldots \otimes \xi_{f_n})  \\
&\quad- \sum_{i=1}^n q^{i-1} \ip{\xi_{f_0}}{\xi_{f_i}} W(\xi_{f_1} \otimes \ldots \otimes \check{\xi}_{f_i} \otimes \ldots \otimes \xi_{f_n}) \\
&\quad- \sum_{i=1}^n q^{i-1} W(T_{f_0} (\xi_{f_i}) \otimes \xi_{f_1} \otimes \ldots \otimes \check{\xi}_{f_i} \otimes \ldots \otimes \xi_{f_n})
\end{split}
\end{equation}
and
\begin{equation}
\label{Wick}
W(\eta_1 \otimes \ldots \otimes \eta_n) \Omega = \eta_1 \otimes \ldots \otimes \eta_n.
\end{equation}

\begin{Prop}
\label{Prop:Separating}
The vacuum vector $\Omega$ is separating for $\Aalg_{X, \Delta}$.
\end{Prop}

\begin{proof}
It is clear from the definition that $W(\eta_1 \otimes \ldots \otimes \eta_n)$ is a polynomial in the $X$ operators (and so lies in $\Aalg_{X, \Delta}$), and conversely that the mapping $W$ is onto. In fact,
\[
W(\xi_{f_1} \otimes \xi_{f_2} \otimes \ldots \otimes \xi_{f_n})
= X(f_1) X(f_2) \ldots X(f_n) + Q,
\]
where $Q$ is a polynomial in the $X$ operators of degree at most $n-1$. Inverting this relation,
\[
X(f_1) X(f_2) \ldots X(f_n)
= W(\xi_{f_1} \otimes \xi_{f_2} \otimes \ldots \otimes \xi_{f_n})
+ W(\vec{\eta}),
\]
where $\vec{\eta} \in \bigoplus_{i=0}^{n-1} \mc{H}^{\otimes i}$. Now let $A \in \Aalg_{X, \Delta}$. Then $A = W(\vec{\eta})$ for some $\vec{\eta} \in \Falg(\mc{D})$. If $A \Omega = 0$, then $W(\vec{\eta}) \Omega = \vec{\eta} = 0$, and so $0 = W(\vec{\eta}) = A$.
\end{proof}

\begin{Remark}
\label{Remark:Commutant}
The following (possibly unbounded) operators commute with $\Aalg_{X, \Delta}$ on the dense domain $\Falg(\mc{D})$: for $f \in \mf{R}_0 \langle \mb{x}, \mc{S} \rangle$, define the operator $X^r(f)$ by
\[
X^r(f) (\eta_1 \otimes \ldots \otimes \eta_n) 
= W(\eta_1 \otimes \ldots \otimes \eta_n) X(f) (\Omega)
= W(\eta_1 \otimes \ldots \otimes \eta_n) \xi_{f}.
\]
Explicitly, their values on the tensors of low order are
\begin{align*}
X^r(f) (\eta) 
&= \xi_f \otimes \eta + T_f \eta + \ip{\eta}{\xi_f}, \\
X^r(f) (\eta_1 \otimes \eta_2)
&= \eta_1 \otimes \eta_2 \otimes \xi_f 
+ \bigl( q \ip{\eta_1}{\xi_f} \eta_2 + \ip{\eta_2}{\xi_f} \eta_1 \bigr) \\
&\quad+ \bigl( q (T_f \eta_1) \otimes \eta_2 + \eta_1 \otimes (T_f \eta_2) \bigr), \\
X^r(f) (\eta_1 \otimes \eta_2 \otimes \eta_3)
&= \eta_1 \otimes \eta_2 \otimes \eta_3 \otimes \xi_f \\
&\quad+ \bigl( q^2 \ip{\eta_1}{\xi_f} \eta_2 \otimes \eta_3 
+ q \ip{\eta_2}{\xi_f} \eta_1 \otimes \eta_3 
+ \ip{\eta_3}{\xi_f} \eta_1 \otimes \eta_2 \bigr) \\
&\quad+ \bigl( q^2 (T_f \eta_1) \otimes \eta_2 \otimes \eta_3 
+ q \eta_1 \otimes (T_f \eta_2) \otimes \eta_3 
+ \eta_1 \otimes \eta_2 \otimes (T_f \eta_3) \bigr) \\
&\quad+ Q(f)(\eta_1 \otimes \eta_2 \otimes \eta_3),
\end{align*}
where
\begin{multline*}
Q(f)(\xi_{g_1} \otimes \xi_{g_2} \otimes \xi_{g_3})
= q(1-q) \bigl[ \ip{T_{g_1} \xi_{g_3}}{\xi_f} \xi_{g_2}
+ T_{g_1} T_{g_3} \xi_f \otimes \xi_{g_2} \\
- \ip{\xi_{g_1}}{\xi_{g_3}} T_{g_2} \xi_f
- T_{g_1} \xi_{g_3} \otimes  T_{g_2} \xi_f \bigr]
\end{multline*}
\end{Remark}

\begin{Notation}
A set partition $\pi$ of a set $\mc{T}$ is a collection of disjoint subsets of $\mc{T}$ whose union equals $\mc{T}$. Let $\Part(n)$ be the collection of all set partitions of the set $\set{1, 2, \ldots, n}$. Let $\pi \in \Part(n)$ be a set partition, $\pi = \set{B_1, B_2, \ldots, B_{\abs{\pi}}}$. Order the classes according to their first elements, that is, $\min(B_1) < \min(B_2) < \ldots  < \min(B_{\abs{\pi}})$. For $S \subset \pi$, call the pair $(S, \pi)$ an extended partition; $S$ is to be thought of as the collection of classes ``open on the left''. See Figure~\ref{Figure:Extended} for an example. 
\begin{figure}[ht]
  \psfig{figure=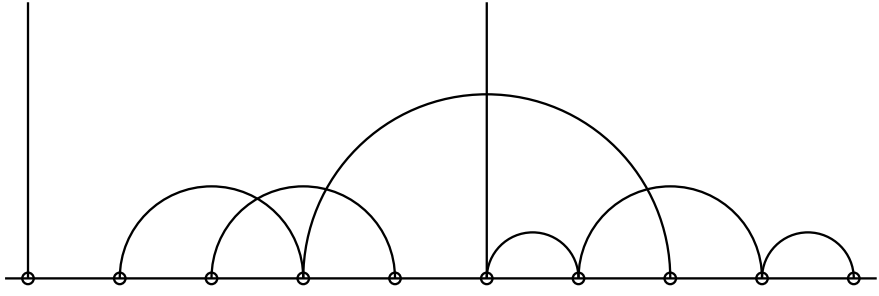,height=1.5in,width=0.7\textwidth}
  \caption{An extended partition of $10$ elements with $2$ left-open classes and $4$ restricted crossings.}
  \label{Figure:Extended}
\end{figure}
For $1 \leq k < m \leq n$, define the restriction
\[
(S', \pi') = (S, \pi) \upharpoonright \set{k, \ldots, m}
\]
as follows:
\begin{align*}
B' \in \pi' &\text{ if } B' = B \cap \set{k, \ldots, m}, B \in \pi, \\
B' \in S' &\text{ if } B \in S \text{ or } B \cap \set{1, \ldots, k-1} \neq \emptyset.
\end{align*}
Define the number of right restricted crossings of $(S, \pi)$ at the point $k$ as follows:
\[
\rc{k, S, \pi} =
\begin{cases}
0, & \text{if } k \in B, k = \max(B), \\
\abs{S'}, & \text{if } k \in B, j = \min \set{i \in B, i > k}, \\
& \qquad (S', \pi') = (S, \pi) \upharpoonright \set{k+1, \ldots, j-1}. 
\end{cases}
\]
Let $\rc{S, \pi} = \sum_{k=1}^n \rc{k, S, \pi}$. Note that also (see \cite{BiaCross})
\[
\rc{S, \pi} 
= \rc{\pi} + \sum_{B \in S} \abs{C \in \pi: \min(C) < \min(B) < \max(C)}.
\]
\end{Notation}

\begin{Prop}
\label{Prop:Product-Wick}
Let $\text{Sing}(\pi)$ denote the single-element classes of a partition $\pi$. Then
\begin{multline}
\label{Wick-product-K}
X(f_1) X(f_2) \ldots X(f_n) \\
= \sum_{\pi \in \Part(n)} \sum_{\text{Sing}(\pi) \subset S \subset \pi} q^{\rc{S, \pi}} \prod_{B \not \in S} \ip{\xi_{f_{\min(B)}}}{\Bigl( \prod_{\substack{i \in B \\ i \neq \min(B), \max(B)}} T_{f_i} \Bigr) \xi_{f_{\max(B)}}} \\
\times W \biggl( \bigotimes_{B \in S} \Bigl( \prod_{\substack{i \in B \\ i \neq \max(B)}} T_{f_i} \Bigr) \xi_{f_{\max(B)}} \biggr).
\end{multline}
\end{Prop}

\begin{proof}
Both sides of the expression~\eqref{Wick-product-K} are in $\Aalg_{X, \Delta}$. Evaluate them on $\Omega$. We obtain
\begin{multline*}
\prod_{i=1}^n \bigl( a(\xi_{f_i}) + a^\ast(\xi_{f_i}) + p(T_{f_i}) \bigr) \Omega \\
= \sum_{\pi \in \Part(n)} \sum_{\text{Sing}(\pi) \subset S \subset \pi} q^{\rc{S, \pi}} \prod_{B \not \in S} \ip{\xi_{f_{\min(B)}}}{\Bigl( \prod_{\substack{i \in B \\ i \neq \min(B), \max(B)}} T_{f_i} \Bigr) \xi_{f_{\max(B)}}} \\ \times \bigotimes_{B \in S} \Bigl( \prod_{\substack{i \in B \\ i \neq \max(B)}} T_{f_i} \Bigr) \xi_{f_{\max(B)}}.
\end{multline*}
For fixed $(S, \pi)$, each term on the right-hand-side, possibly up to a power of $q$, is a summand in the expansion of $Z_1 Z_2 \ldots Z_n \Omega$ for a unique sequence of operators $Z_1, Z_2, \ldots, Z_n$, where
\[
Z_i =
\begin{cases}
a(\xi_{f_i}) & \text{if } i \in B, B \in \pi \backslash S, i = \min(B), \abs{B} > 1, \\
a^\ast(\xi_{f_i}) & \text{if } i \in B, i = \max(B), \text{ either } B \in S \text{ or } (B \in \pi \backslash S, \abs{B} > 1), \\
p(T_{f_i}) & \text{if } i \in B, i \neq \max(B), \text{ either } B \in S \text{ or } (B \in \pi \backslash S, i \neq \min(B)).
\end{cases}
\]
Conversely, any sequence $Z_1 \ldots Z_n$ in the expansion of the left-hand-side of~\eqref{Wick-product-K} whose evaluation on $\Omega$ is non-zero is of this form for some extended partition $(S, \pi)$. It remains to show that its coefficient is exactly $q^{\rc{S, \pi}}$. We show this by induction. Consider a word $Z_1 \ldots Z_n$ corresponding to an extended partition $(S, \pi)$. It is clear that $Z_k Z_{k+1} \ldots Z_n \Omega$ is the element of $\Falg(H)$ corresponding to the restriction $(S, \pi) \upharpoonright \set{k, k+1, \ldots, n}$. By induction hypothesis its coefficient is 
\[
q^{\rc{(S, \pi) \upharpoonright \set{k, k+1, \ldots, n}}}.
\]
If $Z_{k-1} = a^\ast(f_{k-1})$, the coefficient in front of the tensor corresponding to $Z_{k-1} Z_k \ldots Z_n \Omega$ is the same, and so is the number of restricted crossings of $(S, \pi) \upharpoonright \set{k-1, k, \ldots, n}$. On the other hand, suppose $Z_{k-1} = a(f_{k-1})$ or $p(f_{k-1})$. Suppose $k-1 \in B$ and $j = \min \set{i \in B, i > (k-1)}$. Then the degree of $q$ in the coefficient is incremented by
\[
\abs{k-1 < i < j: Z_i = p(f_i) \text{ or } Z_i = a^\ast(f_i)}.
\]
But this is exactly the number of right restricted crossings of $(S, \pi)$ at $(k-1)$.
\end{proof}

\section{Generalized chaos decomposition}
Let $\nu$ be a probability measure on $\mf{R}$ with compact support, and in particular with finite moments
\[
r_{k+2} = \int_{\mf{R}} x^k \,d\nu(x)
\]
(note the shift in the index) of all orders. Then the functional $\psi_\nu[x^k] = r_k$, $\psi_\nu[x] = r_1 = 0$, $\psi_\nu[1] = 0$ on $\mf{C}[x]$ is conditionally positive definite and both left- and right-bounded. Let $\mc{S} = \set{[a, b) \subset \mf{R}_+}$. On the corresponding  $\mf{R} \langle \mb{x}, \mc{S} \rangle$, define the functional
\[
\psi \Bigl[ \prod_{i=1}^n x_{I_i} \Bigr] = \abs{\bigcap_{i=1}^n I_i} r_n.
\]
This functional is also conditionally positive and bounded. The corresponding $\mc{H}$ is naturally isomorphic to $L^2(\mf{R}_+, dt) \otimes V$, where $V \cong \overline{\ls{x^k : k \geq 1}} \subset L^2(\mf{R}, \frac{1}{x^2} \nu(dx))$. Note that the Hilbert space $V$ is isomorphic to $L^2(\mf{R}, \nu)$ via the map $f \mapsto f/x$. So
\[
\mc{H} \cong L^2(\mf{R}_+, dt) \otimes L^2(\mf{R}, \nu).
\]

Denote $X(I) = p(\xi_I, T_I)$. Note that $X(I)$ corresponds to $\chf{I} \otimes 1 \in \mc{H}$. Similarly, $Y_k (I) = p(\xi_{I, I, \ldots, I}, T_{I, I, \ldots, I})$ naturally corresponds to $\chf{I} \otimes x^{k-1}$, where $x$ is the independent variable on $(\mf{R}, \nu)$. Finally, denote $X(t) = X([0, t))$, $Y_k(t) = Y_k([0, t))$. Then all $\set{Y_k(t)}_{t \in [0, \infty)}$ are $q$-L\'{e}vy processes in the sense of \cite{AnsQCum}.

We could have similarly defined multi-dimensional processes, but the notation gets heavier, while the phenomena are the same.

Note that $\Aalg_X$ is equal to the algebra generated by $\set{X(t)}_{t \in [0, \infty)}$, and $\Aalg_{X, \Delta}$ is equal to the algebra generated by $\set{Y_k(t)}_{k \in \mf{N}, t \in [0, \infty)}$. On $\Aalg_X$, define the vacuum state
\[
\state{A} = \ip{\Omega}{A \Omega}
\]
and an inner product
\[
\ip{A}{B}_\phi = \state{A^\ast B} = \ip{A \Omega}{B \Omega}_q.
\]
The completion of $\Aalg_X$ with respect to this inner product will be denoted by $L^2(\Aalg_X, \varphi)$.

\begin{Lemma}
\label{Lemma:Diagonal}
For a subdivision $\mc{I} = \set{I_i}_{i=1}^N$ of $[0, t)$, denote $\delta(\mc{I}) = \max_i \abs{I_i}$. Define the diagonal measure $\Delta_k(t)$ of $X$ to be the limit
\[
\Delta_k(t) = \lim_{\delta(\mc{I}) \rightarrow 0} \sum_{i=1}^N (X(I_i))^k := \int_0^t (dX(t))^k.
\]
Then the limit exists in $L^2(\phi)$, the diagonal measures of $X$ are elements of $L^2(\Aalg_X, \varphi)$, and 
\[
\Delta_k(t) = Y_k(t) + r_k.
\]
In particular, the operators $Y_k(t)$ can in fact be identified with elements of $L^2(\Aalg_X, \varphi)$.
\end{Lemma}

\begin{proof}
See the appendix of~\cite{AnsLinear}.
\end{proof}

Note that $\norm{x^2}_{\psi_\nu} = r_2 = \nu(\mf{R}) = 1$, so the variance of $X$ at time one is $1$.

Since the $q$-inner product is non-degenerate, $L^2(\mf{R}_+, dt)^{\otimes n}$ with the usual and the $q$-inner product are isomorphic as vector spaces. Since $L^2(\mf{R}_+, dt)^{\otimes n} \cong L^2(\mf{R}_+^n, dt^{\otimes n})$, we can induce the $q$-inner product on this space and denote the resulting Hilbert space by $L^2_q(\mf{R}_+^n, dt^{\otimes n})$.

Note also that $\mf{C}\langle x_1, x_2, \ldots, x_n \rangle$ can be naturally identified with a dense subset of $V^{\otimes n}$.

\begin{Defn}
Let $X_1, X_2, \ldots, X_k$ be centered $q$-L\'{e}vy processes normalized to have variance $1$ at time $1$. Let $F$ be an indicator function $F = \chf{I_1 \times I_2 \times \ldots \times I_k}$ such that the intervals $\set{I_i}_{i=1}^k$ are disjoint. Define the multiple Wiener-It\^{o} stochastic integral of $F$
\[
\int_0^\infty F(t_1, t_2, \ldots, t_k) \,dX_1(t_1) \,dX_2(t_2) \ldots \,dX_k(t_k)
\]
to be $\prod_{i=1}^k X_i(I_i)$. Extend the definition to $\mf{C}$-linear combinations of such functions in a linear way. Such an integral will frequently be denoted simply by $\int F \,dX_1 \,dX_2 \ldots \,dX_k$.
\end{Defn}

\begin{Prop}
Let $X_1, X_2, \ldots, X_k$ be as in the preceding definition. The map
\[
F \mapsto \int F \,dX_1 \,dX_2 \ldots \,dX_k
\]
extends to an isometry from $L^2_q(\mf{R}_+^k, dt^{\otimes k})$ into $L^2(\Aalg_{X_1, X_2, \ldots, X_k}, \varphi)$. Here $\Aalg_{X_1, X_2, \ldots, X_k}$ is the algebra generated by $\set{X_i(t)}_{t \in [0, \infty), i  = 1, 2, \ldots, k}$.
\end{Prop}

\begin{proof}
Let $F$ be a simple function, $F = \sum a_{\vec{u}} \chf{I_{u(1)} \times I_{u(2)} \times \ldots \times I_{u(k)}}$, such that the intervals $\set{I_{u(i)}}_{i=1}^k$ are disjoint for each $\vec{u}$. By definition,
\[
\int F \,dX_1 \,dX_2 \ldots \,dX_k = \sum a_{\vec{u}} X_1(I_{u(1)}) X_2(I_{u(2)}) \ldots X_k(I_{u(k)}).
\]
Let $G = \sum b_{\vec{v}} \chf{J_{v(1)} \times J_{v(2)} \times \ldots \times J_{v(k)}}$ be a function of the same type. In the following expression, $\Part_2(k,k)$ denotes the pair partitions $\pi$ of $\set{1, 2, \ldots 2k}$ such that for all $(i, j) = B \in \pi$, $i \leq k, j > k$. Such a partition naturally induces a permutation $\sigma \in Sym(k)$ by $\sigma(i) = j-k$ for $(k+1-i) \stackrel{\pi}{\sim} j$. It is easy to see that $\rc{\pi} = i(\sigma)$.
\begin{align*}
&\ip{\int F \,dX_1 \,dX_2 \ldots \,dX_k}{\int G \,dX_1 \,dX_2 \ldots \,dX_k}_\phi \\
&\quad = \state{\sum \bar{a}_{\vec{u}} b_{\vec{v}} X_k(I_{u(k)}) \ldots X_2(I_{u(2)}) X_1(I_{u(1)}) X_1(J_{v(1)}) X_2(J_{v(2)}) \ldots X_k(J_{v(k)})} \\
&\quad = \sum a_{\vec{u}} b_{\vec{v}} \sum_{\pi \in \Part_2(k,k)} q^{\rc{\pi}} \prod_{i=1}^k \abs{I_{u(k+1-\min(B_i))} \cap J_{v(\max(B_i)-k)}} \\
&\quad = \sum a_{\vec{u}} b_{\vec{v}} \sum_{\sigma \in \Sym(n)} q^{i(\sigma)} \prod_{i=1}^k \abs{I_{u(i)} \cap J_{v(\sigma(i))}} \\
&\quad = \ip{F}{G}_q.
\end{align*}
Since such simple functions are dense in $L^2_q(\mf{R}_+^k, dt^{\otimes k})$, the stochastic integral map can be isometrically extended to the whole of this Hilbert space.
\end{proof}

In particular, define the full stochastic measure
\[
\psi_k((X_1, X_2, \ldots, X_k); t_1, t_2, \ldots, t_k) = \int \chf{\prod_{i=1}^k [0, t_i)} \,dX_1 \,dX_2 \ldots \,dX_k.
\]
If all $X_i$'s are equal to some $X$, we will omit it from the notation, and write simply $\psi_k(t_1, t_2, \ldots, t_k)$. If all $t_i$'s are equal to $t$, we write $\psi_k(t)$.

If the variances of the integrator processes are not normalized to $1$, the multiple stochastic integral map differs from an isometry by a constant factor, and is again well-defined.

\begin{Notation}
Let $X$ be a $q$-L\'{e}vy process. Let $\set{\hat{Y}_k(t)}$ be the Gram-Schmidt orthogonalization of $\set{Y_k(t)}$ in $L^2(\mc{A}_X, \varphi)$. Note that $\state{Y_k(t) Y_j(t)} = t r_{k+j}$ and $r_k = \int_{\mf{R}} x^{k-2} \,d\nu(x)$, where $\nu$ is a ($q$-)canonical measure for this process. Thus the coefficients of
\[
Y_1(t), Y_2(t), \ldots, Y_j(t), \ldots
\]
in the expansion of $\hat{Y_k}(t)$ are precisely those of $1, x, \ldots, x^{j-1}, \ldots$ in the orthogonal polynomials $P_{k-1}$ with respect to the measure $t\, d \nu(x)$. Equivalently, $\hat{Y}_k(I)$ corresponds to $\chf{I} \otimes P_{k-1} \in \mc{H}$. 
\end{Notation}

\begin{Prop}
\label{Prop:Orthogonal}
For a multi-index $\vec{u}$, denote 
\[
H_{\vec{u}} 
= \set{\int F \,d\hat{Y}_{u(1)} \,d\hat{Y}_{u(2)} \ldots \,d\hat{Y}_{u(n)}: F \in L^2_q(\mf{R}_+^n, dt^{\otimes n})}.
\]
Then these subspaces are orthogonal for different $\vec{u}$.
\end{Prop}

\begin{proof}
It suffices to prove the result for simple $F$. In that case the argument is similar to the proof of the preceding proposition. It suffices to note that $\state{\hat{Y}_k(I) \hat{Y}_j(I)} = 0$ for $k \neq j$.
\end{proof}

\begin{Defn}
\label{Defn:St}
Let $\mc{I} = \set{I_i}_{i=1}^N$ be a subdivision of $[0, t)$. For $\pi \in \Part(n)$, define
\[
\St{\pi}(t; \mc{I})
= \sum_{\vec{u} \in [N]^n_\pi} X(I_{u(1)}) X(I_{u(2)}) \ldots X(I_{u(n)}),
\]
where
\[
[N]^n_\pi = \set{\vec{u} \in \set{1, \ldots, N}^n : u(i) = u(j) \Leftrightarrow i \stackrel{\pi}{\sim} j}.
\]
Define the partition-dependent stochastic measures
\[
\St{\pi}(t) = \lim_{\delta(\mc{I}) \rightarrow 0} \St{\pi}(t; \mc{I})
\]
if the limit exists.

Also define $R_\pi(t) = t^{\abs{\pi}} \prod_{B \in \pi} r_{\abs{B}}$. Note that this notation differs from the one in \cite{AnsQCum}.
\end{Defn}

In particular,
\[
\Delta_k(t) = \St{\hat{1}}(t) = \St{\set{(1, 2, \ldots, k)}}(t)
\]
and
\[
\psi_k(t) = \St{\hat{0}}(t) = \St{\set{(1), (2), \ldots, (k)}}(t).
\]

\begin{Prop}
\label{Prop:St_pi}
Partition-dependent stochastic measures of a $q$-L\'{e}vy process are well-defined as limits in $L^2(\Aalg_X, \varphi)$, and equal to
\begin{equation}
\label{St_pi}
\St{\pi}(t)
= \sum_{S \subset \pi} q^{\rc{S, \pi}} R_{\pi \backslash S}(t) \psi(Y_{\abs{B}}(t): B \in S).
\end{equation}
\end{Prop}

\begin{proof}
Let $\mc{I}$ be a subdivision of the interval $[0, t)$. We will show that
\begin{equation}
\label{St_pi;I}
\norm{\St{\pi}(t; \mc{I}) - \sum_{S \subset \pi} q^{\rc{S, \pi}} R_{\pi \backslash S}(t) \psi(Y_{\abs{B}}(t): B \in S; \mc{I})}_2 \rightarrow 0
\end{equation}
as $\delta(\mc{I}) \rightarrow 0$. Since the $L^2$-limit of the second term of \eqref{St_pi;I} exists and equals the right-hand-side of \eqref{St_pi}, this will also show that the limit of the left-hand-side exists. Moreover, since the first term of \eqref{St_pi;I} is in $\Aalg_X$, all the quantities involved are in $L^2(\Aalg_X, \varphi)$

Since $X(I) = a(I) + a^\ast(I) + p(I)$ and the intervals in $\mc{I}$ are disjoint,
\[
\St{\pi}(t; \mc{I}) \Omega 
= \sum_{\vec{u} \in [N]^n_{\pi}} \sum_{\substack{\sigma \in \Part(n) \\ \sigma \leq \pi}} \sum_{S \subset \sigma} a_{(S, \sigma)} \prod_{B \not \in S} \Bigl( \abs{\bigcap_{i \in B} I_{u(i)}} r_{\abs{B}} \Bigr) \bigotimes_{B \in S} \Bigl( \chf{\bigcap_{i \in B} I_{u(i)}} \otimes x^{\abs{B} - 1} \Bigr),
\]
where the sum is over all refinements of $\pi$, and each $a_{(S, \sigma)}$ is a power of $q$. Since $\vec{u} \in [N]^n_{\pi}$ and $\sigma \leq \pi$, we may write $u(B)$ for any $u(i), i \in B$, so the preceding expression is equal to
\begin{equation}
\label{St-I}
\sum_{\vec{u} \in [N]^n_{\pi}} \sum_{\substack{\sigma \in \Part(n) \\ \sigma \leq \pi}} \sum_{S \subset \sigma} a_{(S, \sigma)} \prod_{B \not \in S} \Bigl( \abs{I_{u(B)}} r_{\abs{B}} \Bigr) \bigotimes_{B \in S} \Bigl( \chf{I_{u(B)}} \otimes x^{\abs{B} - 1} \Bigr).
\end{equation}
Also,
\begin{multline}
\label{R-Psi}
\sum_{S \subset \pi} q^{\rc{S, \pi}} R_{\pi \backslash S}(t) \psi(Y_{\abs{B}}(t): B \in S; \mc{I}) \Omega \\
= \sum_{S \subset \pi} q^{\rc{S, \pi}} \prod_{B \not \in S} (r_{\abs{B}} t) \sum_{\vec{v} \in [N]^{\abs{S}}_{\hat{0}}} \bigotimes_{j=1}^{\abs{S}} \Bigl( \chf{I_{v(j)}} \otimes x^{\abs{B_j} - 1} \Bigr).
\end{multline}
As $\delta(\mc{I}) \rightarrow 0$, each term in equation~\eqref{St-I} for $\sigma = \pi$ converges to
\[
a_{(S, \pi)} \prod_{B \not \in S} \Bigl( t \cdot r_{\abs{B}} \Bigr) \bigotimes_{B \in S} \Bigl( \chf{[0, t)} \otimes x^{\abs{B} - 1} \Bigr).
\]
This is also the limit of the corresponding term in equation~\eqref{R-Psi}, provided that we show $a_{(S, \pi)} = q^{\rc{S, \pi}}$. This follows via an argument similar to the one in Proposition~\ref{Prop:Product-Wick}.

On the other hand, the norm of the vector in \eqref{St-I} is bounded by
\[
c \sum_{\vec{u} \in [N]^n_{\pi}} \sum_{\substack{\sigma \in \Part(n) \\ \sigma \leq \pi}} \prod_{B \in \sigma} \abs{I_{u(B)}},
\]
where $c$ is a constant depending on $q, \set{r_i}$ but independent of $\mc{I}, \pi$. For each $\sigma < \pi$, some $C_1, C_2 \in \sigma$ are in the same class of $\pi$. As a result, $u(C_1) = u(C_2)$ for all $\vec{u} \in [N]^n_{\pi}$. We may assume $\delta(\mc{I}) < 1$. Then 
\[
\sum_{\vec{u} \in [N]^n_{\pi}} \prod_{B \in \sigma} \abs{I_{u(B)}}
\leq \sum_{\vec{u} \in [N]^n_{\pi}} \abs{I_{u(C_1)}} \prod_{B \in \pi} \abs{I_{u(B)}}
\leq \delta(\mc{I}) \sum_{\vec{u} \in [N]^n_{\pi}} \prod_{B \in \pi} \abs{I_{u(B)}}
\leq \delta(\mc{I}) t^{\abs{\pi}},
\]
which converges to $0$ as $\delta(\mc{I}) \rightarrow 0$.
\end{proof}

\begin{Ex}
\label{Ex:Classical}
Processes with independent increments correspond to $q=1$. In this case the formula takes the form
\[
\St{\pi}(t)
= \sum_{S \subset \pi} R_{\pi \backslash S}(t) \psi(Y_{\abs{B}}(t): B \in S)
= \psi(\Delta_{\abs{B_1}}(t), \Delta_{\abs{B_2}}(t), \ldots, \Delta_{\abs{B_{\abs{\pi}}}}(t)).
\]
Processes with freely independent increments correspond to $q=0$. In this case the formula takes the form
\[
\St{\pi}(t)
= \sum_{S \subset \text{ Outer}(\pi)} R_{\pi \backslash S}(t) \psi(Y_{\abs{B}}(t): B \in S)
= R_{\text{Inner}(\pi)}(t) \psi(\Delta_{\abs{B}}(t): B \in { \text{Outer}(\pi)}).
\]
for $\pi \in \NC(n)$, and $0$ otherwise. Here $\NC(n)$ are all the non-crossing partitions, and $\text{Inner}(\pi)$, $\text{Outer}(\pi)$ are the inner, respectively, outer classes of $\pi$. This result is the main theorem of \cite{AnsFSM2} (with a weaker mode of convergence).

Finally, in the $q$-Gaussian case, the formula takes the form
\[
\St{\pi}(t)
= \sum_{\substack{S \subset \text{ Sing}(\pi), \\
\pi \backslash S \subset \text{ Pairs}(\pi)}} q^{\rc{S, \pi}} R_{\pi \backslash S}(t) \psi(Y_{\abs{B}}(t): B \in S)
= q^{\rc{\text{Sing}(\pi), \pi}} R_{\text{Pairs}(\pi)}(t) \psi_{\abs{\text{Sing}(\pi)}}(t).
\]
for $\pi \in \Part_{1,2}(n)$, and $0$ otherwise. Here $\text{Pairs}(\pi)$, $\text{Sing}(\pi)$ are $2$, respectively, $1$-element classes of $\pi$. Note that in this case, $\rc{\text{Sing}(\pi), \pi} = \rc{\pi} +$ the singleton depth of $\pi$, $R_{\text{Pairs}(\pi)}(t) = t^{\abs{\text{Pairs}(\pi)}}$. Thus we recover Proposition 6.12 of \cite{AnsQCum} (with a weaker mode of convergence). For closely related results, see \cite{Effros-Popa}.
\end{Ex}

\begin{Ex}
In the centered $q$-Charlier case,
\[
\St{\pi}(t)
= \sum_{\text{Sing}(\pi) \subset S \subset \pi} q^{\rc{S, \pi}} t^{n - \abs{S}} C_{\abs{S}, q}(X(t),t).
\]
As a consequence,
\[
x^n = \sum_{\pi \in \Part(n)} \sum_{\text{Sing}(\pi) \subset S \subset \pi} q^{\rc{S, \pi}} t^{n - \abs{S}} C_{\abs{S}, q}(X(t),t).
\]
Here $C_{n, q}$ are the continuous big $q$-Hermite polynomials, see for example \cite{AnsLinear}.
\end{Ex}

\begin{Lemma}
\label{Lem:Dense}
The image of the multiple stochastic integral map with respect to the processes $\set{Y_k}$ contains $\Aalg_{X, \Delta}$.
\begin{enumerate}
\item
\label{Part:Poly}
Polynomials in $X$ can be expressed as multiple stochastic integrals:
\[
X^n(t) = \sum_{\pi \in \Part(n)} \St{\pi}(t).
\]
\item
\label{Part:Cyclic}
For $F \in L^2_q(\mf{R}_+^k, dt^{\otimes k})$,
\begin{equation*}
W(F \otimes (P_{u(1)-1}(x_1) P_{u(2)-1}(x_2) \ldots P_{u(k)-1}(x_k))) 
= \int F \,d\hat{Y}_{u(1)} \,d\hat{Y}_{u(2)} \ldots \,d\hat{Y}_{u(k)}
\end{equation*}
\end{enumerate}
\end{Lemma}

\begin{proof}
The first part of the lemma is a basic, purely combinatorial, property of partition-dependent stochastic measures due to \cite{Rota}. For the second part, we observe that the maps
\[
F \mapsto \int F \,d\hat{Y}_{u(1)} \,d\hat{Y}_{u(2)} \ldots \,d\hat{Y}_{u(k)}
\]
and $A \mapsto A \Omega$ are isometries. Therefore it suffices to show the property for simple functions
\[
F = \chf{I_1 \times I_2 \times \ldots \times I_k},
\]
where all $I_j$ are disjoint. But in this case
\begin{align*}
\int F \,d\hat{Y}_{u(1)} \,d\hat{Y}_{u(2)} \ldots \,d\hat{Y}_{u(k)}
&= \hat{Y}_{u(1)}(I_1) \hat{Y}_{u(2)}(I_2) \ldots \hat{Y}_{u(k)}(I_k) \\
&= W(F \otimes (P_{u(1)-1}(x_1) P_{u(2)-1}(x_2) \ldots P_{u(k)-1}(x_k))).
\end{align*}
\end{proof}

\begin{Cor}
\label{Cor:Wick-integrals}
In particular,
\[
W \bigl(\chf{\prod_{i=1}^k [0, t_i)} \otimes (x_1^{u(1)} x_2^{u(2)} \ldots x_k^{u(k)}) \bigr) = \psi((Y_{u(1)+1}, Y_{u(2)+1}, \ldots, Y_{u(k)+1}); t_1, t_2, \ldots, t_k).
\]
\end{Cor}

\begin{Cor}
The vacuum vector $\Omega$ is cyclic and separating for $\Aalg_X$. The representation $A \mapsto A \Omega$ of $\Aalg_X$ is faithful, and extends to an isomorphism between $L^2(\Aalg_X, \varphi)$ and $\mc{F}_q(\mc{H})$. The Wick map extends to the inverse isomorphism.
\end{Cor}

\begin{proof}
We only need to show that $\Omega$ is cyclic for $\Aalg_X$. By Proposition~\ref{Prop:St_pi}, for a fixed $\pi$
\[
\psi(Y_{\abs{B}}(t): B \in \pi)
= \St{\pi}(t) - \sum_{\substack{S \subset \pi \\ S \neq \pi}} q^{\rc{S, \pi}} R_{\pi \backslash S}(t) \psi(Y_{\abs{C}}(t): C \in S).
\]
So by induction, 
\[
\psi((Y_{u(1)}, Y_{u(2)}, \ldots, Y_{u(k)}); t_1, t_2, \ldots, t_k) \in L^2(\Aalg_X, \varphi).
\]
Therefore, using the same proposition, each
\[
\St{\pi}((Y_{u(1)}, Y_{u(2)}, \ldots, Y_{u(k)}); t_1, t_2, \ldots, t_k) \in L^2(\Aalg_X, \varphi),
\]
where this expression is defined in the obvious way. Since, similarly to Lemma~\ref{Lem:Dense},
\[
\prod_{k=1}^n Y_{u(k)}(t_k) = \sum_{\pi \in \Part(n)} \St{\pi}((Y_{u(1)}, Y_{u(2)}, \ldots, Y_{u(n)}); t_1, t_2, \ldots, t_n),
\]
we conclude that $\Aalg_{X, \Delta} \subset L^2(\Aalg_X, \varphi)$. But $\Omega$ is cyclic for this algebra.
\end{proof}

The following proposition is an analog of a result of \cite{SchChaotic}.

\begin{Prop}
Any $A \in L^2(\Aalg_X, \varphi)$ has a unique chaos decomposition
\[
A = \sum_{n=0}^\infty \sum_{\vec{u}} \int F_{\vec{u}} \,d\hat{Y}_{u(1)} \,d\hat{Y}_{u(2)} \ldots \,d\hat{Y}_{u(n)},
\]
where
\[
\norm{A}_2^2 = \sum_{\vec{u}} \norm{F_{\vec{u}}}_2^2
\]
and
\[
F_{\vec{u}} \in L^2_q(\mf{R}^n_+, dt^{\otimes n}).
\]
Conversely, any such series converges to an element of $L^2(\Aalg_X, \varphi)$.
\end{Prop}

\begin{proof}
By Proposition~\ref{Prop:Orthogonal}, it suffices to show that $\bigoplus_{\vec{u}} H_{\vec{u}}$ is dense in $L^2(\Aalg_X, \varphi)$. Using the isomorphism between $L^2(\Aalg_X, \varphi)$ and $\mc{F}_q(\mc{H})$, this is guaranteed by Lemma~\ref{Lem:Dense}~\eqref{Part:Cyclic}.
\end{proof}

\begin{Remark}[Classical L\'{e}vy processes]
\label{Remark:Classical}
A few modifications are necessary for the classical case $q=1$. From the point of view of $q$-deformations of the full Fock space, this case is degenerate, since it involves the reduction to the symmetric Fock space. An easy way to modify the preceding arguments for this context is to work, instead of functions in $L^2(\mf{R}^n_+, dt^{\otimes n})$, with square-integrable functions with support in the simplex
\[
D_n = \set{t_1 > t_2 > t_3 > \ldots > t_n \geq 0}.
\]
In this case the von Neumann algebras are commutative, and as a result instead of working with bounded self-adjoint operators we can work with essentially self-adjoint operators with an invariant dense domain consisting of analytic vectors. For the measure $\nu$, this corresponds to dropping the requirement of compact support and instead requiring that it has a finite moment generating function,
\[
\int_{\mf{R}} e^{\theta x} \,d\nu(x) < \infty
\]
for $\theta$ small enough. Note that this is exactly the hypothesis of \cite{SchChaotic}. Moreover, the vacuum vector $\Omega$ is cyclic and separating for the von Neumann algebra (and not just the algebra) generated by $\set{X(t)}$. With these modification, all the preceding statements about isometries and orthogonality remain true, and the result of \cite{SchChaotic} follows. In this case all the operators $\set{X(t)}$ commute and are independent with respect to the expectation $\phi$, and the corresponding convolution semigroup is given by
\[
\log \mc{F}(\mu_t)(\theta) = t \int_{\mf{R}} (e^{i \theta x} - 1 - i \theta x) \frac{1}{x^2} d \nu(x).
\]
Note also that by combining Lemma~\ref{Lem:Dense}~\eqref{Part:Poly} with Proposition~\ref{Prop:St_pi} (more specifically, with Example~\ref{Ex:Classical}), we get an explicit formula
\[
\prod_{i=1}^n X([a_i, b_i)) = \sum_{\pi \in \Part(n)} \int_{\max_{i \in B_1} a_i}^{\min_{i \in B_1} b_i} \cdots \int_{\max_{i \in B_{\abs{\pi}}} a_i}^{\min_{i \in B_{\abs{\pi}}} b_i} d\Delta_{B_1}(t_1) \,d\Delta_{B_2}(t_2) \ldots \,d\Delta_{B_{\abs{\pi}}}(t_{\abs{\pi}}).
\]
It is to be compared with the results of Section 3.1 of~\cite{SchChaotic}. Note also a number of previous results in this direction, such as \cite{Kailath-Segall,ItoSpectral,SegalTensorAlgebras}.
\end{Remark}

\begin{Remark}[Free L\'{e}vy processes]
\label{Remark:Free}
In the free case $q=0$, the von Neumann algebras are no longer commutative, but the vacuum expectation is tracial (which is not the case in general: see the next section). As a result, we can again work with unbounded operators, and weaken the hypothesis to a measure with a finite moment-generating function. Using the operators from Remark~\ref{Remark:Commutant}, we can again show that the vacuum vector is separating for the von Neumann algebra of $\set{X(t)}$. The processes have freely independent increments, and the corresponding free convolution semigroup is given by
\[
z R_{\mu_t}(z) = t \int_{\mf{R}} \left( \frac{1}{1 - z x} - 1  - z x \right) \frac{1}{x^2} \,d\nu(x),
\]
where $R_{\mu}$ is the $R$-transform of $\mu$. 
\end{Remark}

\begin{Cor}
\label{Cor:Step1}
Let $\pi = \set{(1, u(2), u(3), \ldots, u(k-1), n), (2), (3), \ldots, (n-1)}$. Then
\[
\St{\pi}(t) = q^{n-k} \psi(\Delta_k(t), X(t), \ldots, X(t)).
\]
\end{Cor}

\subsection{Kailath-Segall polynomials}
\begin{Defn}
Let $\set{x_i}_{i=1}^\infty$ be (possibly non-commuting) indeterminates. For $\vec{u} \in \mf{N}^\infty$, define the polynomial $A_{\vec{u}}$ of total degree $\abs{u}$ in the variables
\[
\set{x_j: j = \sum_{i \in S} u(i), S \subset \set{1, \ldots, \abs{u}}}
\]
by the recursion
\begin{equation}
\label{Appell}
A_{(j, \vec{u})}
= x_{j} A_{\vec{u}} 
- \sum_{i=1}^n q^{i-1} r_{j + u(i)} A_{\vec{u} \backslash u(i)} - \sum_{i=1}^n q^{i-1} A_{(j + u(i), \vec{u} \backslash u(i))}
\end{equation}
with initial conditions $A_{\emptyset} = 1$, $A_i = x_i$.
\end{Defn}

These polynomials have apparently not been considered explicitly before; they are in some weak sense analogs of the  Appell polynomials. Because of \cite{Kailath-Segall}, it is appropriate to call them Kailath-Segall polynomials.

From equation~\eqref{Wick-recursion} and Corollary~\ref{Cor:Wick-integrals},
\begin{multline*}
A_{\vec{u}} \left( x_j = Y_j \right)
= W \bigl( (\chf{[0, 1)} \otimes x^{u(1)-1}) \otimes \ldots \otimes (\chf{[0, t)} \otimes x^{u(n)-1}) \bigr) \\
= \int_{[0,t)^n} dY_{u(1)}(t_1) dY_{u(2)}(t_2) \ldots dY_{u(n)}(t_n)
\end{multline*}
and similarly,
\[
\int_{\prod_{i=1}^n [0, v(i))} dY_{u(1)}(t_1) dY_{u(2)}(t_2) \ldots dY_{u(n)}(t_n)
\]
are polynomials in $\set{Y_j(\min_{i \in S} v(i)): j = \sum_{i \in S} u(i), S \subset \set{1, \ldots, \abs{u}}}$.

The following proposition is closely related to the results of \cite{Kailath-Segall} in the classical case. We use the notation $[0]_q = 0$, $[n]_q = \sum_{i=0}^{n-1} q^k$,  $[0]_q! = 1$, $[n]_q! = \prod_{i=1}^n [i]_q$.

\begin{Prop}[Kailath-Segall formulas].

\begin{enumerate}
\item
Denote $A^{(n)} = A_{(1,1,\ldots, 1)}$. For a polynomial of total degree $n+1$,
\[
A_{j,1,\ldots,1} = x_j A^{(n)} + \sum_{k=1}^n (-1)^k \frac{[n]_q!}{[n-k]_q!} (x_{j+k} + r_{j+k}) A^{(n-k)}.
\]
\item
Since $r_1 = 0$, 
\[
A^{(n+1)} = \sum_{k=0}^n (-1)^k \frac{[n]_q!}{[n-k]_q!} (x_{k+1} + r_{k+1}) A^{(n-k)}.
\]
\item
In particular,
\[
\psi_n(t)
= \sum_{k=0}^n (-1)^k \frac{[n]_q!}{[n-k]_q!} \Delta_{k+1}(t) \psi_{n-k}(t).
\]
\end{enumerate}
\end{Prop}

\begin{proof}
For the first part, denote the polynomial $A_{j,1,\ldots,1}$ of degree $n+1$ by $A_j^{(n)}$. The result follows by induction from the recursion \eqref{Appell}: $A_j = x_j + r_j$, and
\begin{align*}
A_j^{(n)} 
&= x_{j} A^{(n)} - [n]_q r_{j+1} A^{(n-1)} - [n]_q A_{j+1}^{(n-1)} \\
&= x_{j} A^{(n)} - [n]_q r_{j+1} A^{(n-1)} \\
&\quad - [n]_q \left( x_{j+1} A^{(n-1)} + \sum_{k=0}^{n-1} (-1)^k \frac{[n-1]_q!}{[n-1-k]_q!} (x_{j+k+1} + r_{j+k+1}) A^{(n-k-1)} \right) \\
&= x_{j} A^{(n)} + \sum_{k=1}^n (-1)^k \frac{[n]_q!}{[n-k]_q!} (x_{j+k} + r_{j+k}) A^{(n-k)}.
\end{align*}
The second part follows from the first one, and the third one is an application of Corollary~\ref{Cor:Wick-integrals}. It can also be obtained using Corollary~\ref{Cor:Step1}.
\end{proof}

\begin{Ex}
In the $q$-Gaussian case, $r_2 = 1$, $r_k = 0$ for $k > 2$. Let $x_1 = x$, $x_2 = 1$, $x_k = 0$ for $k > 2$. Then
\[
A^{(n+1)}(x) = x A^{(n)}(x) - [n]_q A^{(n-1)}(x).
\]
So $A^{(n)}$ are the continuous (Rogers) $q$-Hermite polynomials.

In the $q$-Poisson case, $r_k = 1$ for $k \geq 2$. Let $x_k = x$ for $k \geq 1$. Then
\[
A^{(n+1)}(x) = x A^{(n)}(x) - [n]_q A^{(n-1)}(x) - [n]_q A^{(n)}(x).
\]
So $A^{(n)}$ are the centered continuous big $q$-Hermite polynomials, which in this context are $q$-analogs of the Charlier polynomials.
\end{Ex}

See \cite{AnsAppell} for further results in this direction.

\section{Von Neumann algebra of a $q$-L\'{e}vy process}
Let $\mc{A}_{X, \Delta}(t)$ be the von Neumann algebra generated by $\set{Y_k(s): k \in \mf{N}, s < t}$. Let $P_t$ be the orthogonal projection from $L^2(\mf{R}_+, dt)$ onto $L^2([0, t), dt)$. Let
\[
\mc{H}(t) = L^2([0, t), dt) \otimes V.
\]
Denote also by $P_t$ the induced projections from $\mc{H}$ to $\mc{H}(t)$ and from $\Falg(\mc{H})$ to $\Falg(\mc{H}(t))$. It extends to an orthogonal projection from $\mc{F}_q(\mc{H})$ to $\mc{F}_q(\mc{H}(t))$. Note that $\mc{A}_{X, \Delta}(t)$ is also generated as a von Neumann algebra by $W(\Falg(\mc{H}(t)))$.

Denote $E_t[A]= P_t A P_t$. Then $E_t$ is a norm- and strongly continuous projection. Since is clearly maps $W(\Falg(\mc{H}))$ onto $W(\Falg(\mc{H}(t)))$, it maps $\Aalg_{X, \Delta}$ onto $\Aalg_{X, \Delta}(t)$ and $\mc{A}_{X, \Delta}$ onto $\mc{A}_{X, \Delta}(t)$. As a result, $E_t$ is a conditional expectation from $\mc{A}_{X, \Delta}$ to $\mc{A}_{X, \Delta}(t)$, which preserves the vacuum state. In particular, for $A_t, B_t \in \mc{A}_{X, \Delta}(t)$, $Z \in \mc{A}_{X, \Delta}$, $E_t[A_t Z B_t] = A_t E_t[Z] B_t$.

A map $U: \mf{R}_+ \rightarrow \mc{A}_{X, \Delta}$ will be called a process. A simple process is piecewise constant and zero at infinity, $U(t) = \sum_{i=1}^n U_i \chf{[a_i, b_i)}(t)$, where without loss of generality  all the intervals $[a_i, b_i)$ are disjoint. A simple adapted process is a simple process with each $U_i \in \mc{A}_{X, \Delta}(a_i)$. Denote by $\mc{B}_2$ the completion of the set of simple adapted processes with respect to the norm coming from the inner product
\[
\ip{U}{V} = \int_0^\infty \ip{U(t)}{V(t)}_{\phi} dt,
\]
and call the elements of this completion adapted processes. An simple algebraic process is $U(t) = \sum_{i=1}^n U_i \chf{[a_i, b_i)}(t)$ with all $U_i \in \Aalg_{X, \Delta}$, and simple adapted algebraic processes and adapted algebraic processes are defined similarly.

\begin{Defn}
Let $X$ be a $q$-L\'{e}vy process and $U$ a simple adapted process, $U(t) = \sum_{i=1}^n U_i \chf{[a_i, b_i)}(t)$. Define the left and right It\^{o} stochastic integrals
\[
\int_0^\infty U(t) dX(t) = \sum_{i=1}^n U_i X([a_i, b_i))
\]
and
\[
\int_0^\infty dX(t) U(t) = \sum_{i=1}^n X([a_i, b_i)) U_i.
\]
\end{Defn}

Note that for all the arguments below, it is not necessary that the $q$-L\'{e}vy process $\set{X(t)}$ generate the filtration $\set{\mc{A}_{X, \Delta}(t)}$, but only that it be a martingale with respect to it, in other words that $X(I) \Omega \perp \mc{H}(t)$ if $I \cap [0, t) = \emptyset$. In particular, for processes adapted with respect to $\mc{A}_{X, \Delta}(t)$, all the integrals with respect to $\set{Y_k(t)}$, $\set{\Delta_k(t)}$ are defined.

\begin{Lemma}
The stochastic integral map in the preceding definition is $r_2$ times an isometry from $\mc{B}_2$ to $L^2(\mc{A}_{X, \Delta}, \phi)$. Therefore it can be extended to all adapted processes.
\end{Lemma}

\begin{proof}
We consider the right integrals. For two simple adapted processes, 
\begin{equation*}
\ip{\sum_{i=1}^n X([a_i, b_i)) U_i \Omega}{\sum_{j=1}^m X([c_j, d_j)) V_j \Omega}_q 
= \sum_{i=1}^n \sum_{j=1}^m \ip{X([a_i, b_i)) U_i \Omega}{X([c_j, d_j)) V_j \Omega}_q.
\end{equation*}
Since $U, V$ are adapted,
\[
\ip{X([a_i, b_i)) U_i \Omega}{X([c_j, d_j)) V_j \Omega}_q
= \ip{(\chf{[a_i, b_i)} \otimes x) \otimes (U_i \Omega)}{(\chf{[c_j, d_j)} \otimes x) \otimes (V_j \Omega)}_q.
\]
If, say, $b_i < c_j$, this is $0$. On the other hand, if $a_i = c_j, b_i = d_j$, then the inner product is
\[
\abs{b_i - a_i} r_2 \ip{U_i \Omega}{V_j \Omega}_q
= \abs{b_i - a_i} r_2 \ip{U_i}{V_j}_{\phi}.
\]
So for two simple adapted processes,
\[
\ip{\int_0^\infty dX(t) U(t)}{\int_0^\infty dX(t) V(t)}_{\phi} = r_2 \int_0^\infty \ip{U(t)}{V(t)}_{\phi} \,dt.
\]
The result follows. For the left integrals, the proof is similar.
\end{proof}

Denote by $\Gamma_q(q)$ the unitary second quantization of the operator $q \id$, determined by
\[
\Gamma_q(q)(W(\eta_1 \otimes \eta_2 \otimes \ldots \otimes \eta_n)) = q^n W(\eta_1 \otimes \eta_2 \otimes \ldots \otimes \eta_n).
\]
It is a completely positive contraction. The role of this operator has been emphasized by Donati-Martin in \cite{Donati-Martin}.

\begin{Prop}
Let $U \in \Aalg_{X, \Delta}$. Define
\[
\int_0^\infty d X(t) \chf{[u, v)}(t) U d X(t) = \lim_{\delta(\mc{I}) \rightarrow 0} \sum_{i=1}^n X(I_i) U X(I_i),
\]
where $\mc{I} = \set{I_i}$ is a subdivision of $[u, v)$. This limit exists in $L^2(\phi)$. The definition extends linearly to simple adapted algebraic processes so that
\begin{equation}
\label{XUX}
\int_0^\infty d X(t) U(t) d X(t) = \int_0^\infty d \Delta_2(t) \Gamma_q(q)(U(t)).
\end{equation}
\end{Prop}

\begin{proof}
First we prove formula~\eqref{XUX} for an elementary process. Let $U$ be simple, $U = W(\eta_1 \otimes \eta_2 \otimes \ldots \otimes \eta_k) \chf{[u,v)}$, with all the vectors $\eta_i \in \mc{H}(u)$. Let $\mc{I}$ be a subdivision of $[u,v)$. Then by definition, 
\begin{equation}
\label{XWX}
\int_0^\infty d X(t) U(t) d X(t)
= \lim_{\delta(\mc{I}) \rightarrow 0} \sum_{i=1}^N X(I_i) W(\eta_1 \otimes \eta_2 \otimes \ldots \otimes \eta_k) X(I_i).
\end{equation}
Note that for $\zeta \perp \set{\eta_1, \eta_2, \ldots, \eta_k}$, 
\[
W(\eta_1 \otimes \eta_2 \otimes \ldots \otimes \eta_k) \zeta = \eta_1 \otimes \eta_2 \otimes \ldots \otimes \eta_k \otimes \zeta.
\]
Therefore representing the right-hand-side of equation~\eqref{XWX} on $\Omega$, we obtain
\begin{align*}
&\sum_{i=1}^N X(I_i) W(\eta_1 \otimes \eta_2 \otimes \ldots \otimes \eta_k) X(I_i) \Omega \\
&\quad= \sum_{i=1}^N X(I_i) (\eta_1 \otimes \eta_2 \otimes \ldots \otimes \eta_k \otimes (\chf{I_i} \otimes 1)) \\
&\quad= \sum_{i=1}^N \Bigl[ (\chf{I_i} \otimes 1) \otimes \eta_1 \otimes \eta_2 \otimes \ldots \otimes \eta_k \otimes (\chf{I_i} \otimes 1) \\
&\quad\quad + q^k \abs{I_i} r_2 \eta_1 \otimes \eta_2 \otimes \ldots \otimes \eta_k + q^k (\chf{I_i} \otimes x) \otimes \eta_1 \otimes \eta_2 \otimes \ldots \otimes \eta_k \Bigr] \\
&\quad= \sum_{i=1}^N (\chf{I_i} \otimes 1) \otimes \eta_1 \otimes \eta_2 \otimes \ldots \otimes \eta_k \otimes (\chf{I_i} \otimes 1) \\
&\quad\quad+ (v-u) r_2 q^k \eta_1 \otimes \eta_2 \otimes \ldots \otimes \eta_k + q^k (\chf{[u,v)} \otimes x) \otimes \eta_1 \otimes \eta_2 \otimes \ldots \otimes \eta_k ].
\end{align*}
Since the sum in the last term converges to $0$, it follows that
\[
\begin{split}
\int_0^\infty d X(t) U(t) d X(t)
& = q^k \int_0^\infty d Y_2(t) U(t) + q^k r_2 \int_0^\infty U(t) \,dt \\
& = q^k \int_0^\infty d \Delta_2(t) U(t)
= \int_0^\infty d \Delta_2(t) \Gamma_q(q)(U(t)).
\end{split}
\]
By linearity, the same result holds for simple algebraic processes. 
\end{proof}

\begin{Prop}
The algebra $\mc{A}_{X, \Delta}$ has no normal tracial states except, possibly, the vacuum state.
\end{Prop}

\begin{proof}
Suppose $\tau$ is a tracial state on $\mc{A}_{X, \Delta}$. By the previous proposition, for $\set{\eta_j} \in \mc{H}(u)$
\[
\tau[\int_u^v d X(s) W(\eta_1 \otimes \eta_2 \otimes \ldots \otimes \eta_k) d X(s)]
= q^k \tau[\Delta_2([u,v)) W(\eta_1 \otimes \eta_2 \otimes \ldots \otimes \eta_k)].
\]
$\tau$ is normal, so strongly continuous, so $L^2(\phi)$ continuous. Using this and the trace property, the expression above is also equal to
\[
\tau[\int_u^v d X(s) W(\eta_1 \otimes \eta_2 \otimes \ldots \otimes \eta_k) d X(s)]
= \tau[\Delta_2([u,v)) W(\eta_1 \otimes \eta_2 \otimes \ldots \otimes \eta_k)].
\]
Therefore $\tau[\Delta_2([u,v)) W(\eta_1 \otimes \eta_2 \otimes \ldots \otimes \eta_k)] = 0$ for $k > 0$. But note that the limit
\[
\lim_{v \rightarrow \infty} \frac{1}{v-u} \Delta_2([u,v)) = r_2 \id.
\]
Thus $\tau[W(\eta_1 \otimes \eta_2 \otimes \ldots \otimes \eta_k)] = 0$ for $k > 0$.
\end{proof}

\begin{Prop}
The following dichotomy holds.
\begin{enumerate}
\item
For the $q$-Brownian motion, or for $q=0$, the von Neumann algebra is a $\mathrm{II}_1$-factor.
\item
\label{Part:non-BM}
For $q \neq 0$ and all other $q$-L\'{e}vy processes, the von Neumann algebra has not normal tracial states, and thus is an infinite algebra.
\end{enumerate}
\end{Prop}

\begin{proof}
Only part \eqref{Part:non-BM} is new. Let $I \cap J = \emptyset$. Note that 
\[
\state{X(I) X(J) X(I) X(J) Y_k(I)}
= q^2 r_2 r_{2+k} \abs{I} \cdot \abs{J},
\]
while
\[
\state{Y_k(I) X(I) X(J) X(I) X(J)} 
= q r_2 r_{2+k} \abs{I} \cdot \abs{J}.
\]
Thus if the vacuum state is tracial and $q \neq 0$, $r_k = 0$ for $k > 2$. This characterizes the $q$-Brownian motion.
\end{proof}

It is not known if these von Neumann algebras are always factors.

\subsection{Integrals of bi-processes}
Using an idea of Donati-Martin, we can define two-sided stochastic integrals. For $\set{A^i, B^i}_{i=1}^n$ simple algebraic processes, a simple algebraic bi-process is $U = \sum_{i=1}^n A^i \otimes B^i$. $U$ is adapted if each of $A^i, B^i$ is. Denote
\[
\ll A_1 \otimes B_1, A_2 \otimes B_2 \gg = \state{B_1^\ast \Gamma_q(q)(A_1^\ast A_2) B_2}.
\]
This is a positive sesquilinear form (cf.~\cite{Donati-Martin}, although our notation is slightly different). Denote the closure of the space of all adapted simple algebraic bi-processes with respect to the corresponding (semi)norm by $\mc{P}_q$, and call its elements adapted algebraic bi-processes.

\begin{Defn}
For a simple algebraic adapted bi-process $U = \sum_{i=1}^n A^i \otimes B^i$, with $A^i = \sum_{j=1}^N A_j^i \chf{I_j}$, $B^i = \sum_{j=1}^N B_j^i \chf{I_j}$, define the stochastic integral of $U$ to be the operator
\[
\int_0^\infty U(t) \sharp \,dX(t) = \sum_{i=1}^n \sum_{j=1}^N A_j^i X(I_j) B_j^i.
\]
\end{Defn}

\begin{Prop}
The stochastic integral map is an isometry from $\mc{P}_q$ to $L^2(\Aalg_{X, \Delta}, \phi)$.
\end{Prop}

The proof is the same as in the $q$-Brownian motion case, and relies on the key

\begin{Lemma}
For $s < t$ and $Z \in \Aalg_{X, \Delta}(s)$, 
\[
\E_s[X([s,t)) Z X([s,t))] = (t-s) \Gamma_q(q)(Z).
\]
\end{Lemma}

\begin{proof}
Let $Z = W(\eta_1 \otimes \eta_2 \otimes \ldots \otimes \eta_k)$, with $\set{\eta_i} \subset \mc{H}(s)$. Then
\[
\E_s[X([s,t)) Z X([s,t))]
= P_s X([s,t)) Z X([s,t)) P_s
= P_s a([s,t)) Z a^\ast([s,t)) P_s.
\]
Evaluating this expression on $\Omega$, we get
\begin{equation*}
\begin{split}
P_s a([s,t)) (Z \Omega \otimes a^{\ast}([s,t)) \Omega)
&= P_s a([s,t)) (\eta_1 \otimes \eta_2 \otimes \ldots \otimes \eta_k \otimes a^{\ast}([s,t)) \Omega) \\
&= q^k (t-s) P_s (\eta_1 \otimes \eta_2 \otimes \ldots \otimes \eta_k)
= q^k (t-s) P_s Z \Omega.
\end{split}
\end{equation*}
Since the vacuum vector is separating for $\Aalg_{X, \Delta}$ and the conditional expectation maps this algebra into itself, the result follows.
\end{proof}

\appendix

\section{Wick products in the $q$-compound Poisson algebra}
Let $\mc{A}$ be a unital star-algebra, given in a faithful star-representation (by bounded operators) on a Hilbert space $H$ with a cyclic, separating vector $\Omega$. Denote by $\Exp{f} = \ip{\Omega}{f \Omega}$ the state on $\mc{A}$. Construct $\Falg(H)$ and $\mc{F}_q(H)$ as in the beginning of Section~\ref{Sec:q-Fock}, and let $\state{X} = \ip{\Omega}{X \Omega}$ be the vacuum expectation on $\mc{B}(\mc{F}_q(H))$. Note that now, unlike in the body of the paper, we identify $\Omega \in H$ with the vacuum vector in $\mc{F}_q(H)$. We will also identify each $f \in \mc{A}$ with the corresponding vector $f \Omega \in H$.

For $f \in \mc{A}^{sa}$ a self-adjoint element, let
\begin{align*}
& a^\ast(f) \Omega = f, \\
& a^\ast(f) (g_1 \otimes \ldots \otimes g_n) = f \otimes g_1 \otimes \ldots \otimes g_n, \\
& a(f) \Omega = 0, \\
& a(f) g = \Exp{f g} \Omega, \\
& a(f) (g_1 \otimes \ldots \otimes g_n) = \sum_{k=1}^n q^{k-1} \Exp{f g_k} g_1 \otimes \ldots \otimes \check{g}_k \otimes \ldots \otimes g_n, \\
& p(f) \Omega = 0, \\
& p(f) (g_1 \otimes \ldots \otimes g_n) = \sum_{k=1}^n q^{k-1} (f g_k) \otimes g_1 \otimes \ldots \otimes \check{g}_k \otimes \ldots \otimes g_n
\end{align*}
be the creation, annihilation, and gauge operators. Then
\[
X(f) = a(f) + a^\ast(f) + p(f) + \Exp{f}
\]
is a bounded self-adjoint operator.

For $q=1$, if $\mu$ is the distribution of $X(f)$ with respect to $\phi$ and $\nu$ is the distribution of $f$ with respect to $\Exp{\cdot}$, then
\[
\log \int_{\mf{R}} e^{i \theta x} \,d\mu(x) = \int_{\mf{R}} (e^{i \theta x} - 1) \,d\nu(x).
\]
Thus $X(f)$ has a compound Poisson distribution. In particular, if $f$ is a projection, $X(f)$ has a Poisson distribution.

Similarly, for $q=0$ (see \cite{BNTLaws}),
\[
z R_{\mu}(z) = \int_{\mf{R}} \left( \frac{1}{1 - z x} - 1 \right) \,d\nu(x),
\]
where $R_\mu$ is the $R$-transform of $\mu$. So $X(f)$ has a free compound Poisson distribution.

For general $q$ (see \cite{AnsQCum}),
\begin{equation}
\label{Moment}
\state{X(f)^n} = \sum_{\pi \in \Part(n)} q^{\rc{\pi}} \prod_{B \in \pi} \Exp{f^{\abs{B}}}. 
\end{equation}

Let $\Gamma(\mc{A})$ be the algebra generated by $\set{X(f): f \in \mc{A}^{sa}}$, with the obvious involution. As in Lemma~\ref{Lemma:Cyclic} and Proposition~\ref{Prop:Separating}, $\Omega$ is a cyclic and separating vector for $\Gamma(\mc{A})$, so the vacuum state is faithful. Define the Wick product 
\[
W(f_1 \otimes f_2 \otimes \ldots \otimes f_n) \Omega = f_1 \otimes f_2 \otimes \ldots \otimes f_n.
\]
Then $W(f) = X(f) - \Exp{f}$,
\begin{equation}
\label{Wick-recursion2}
\begin{split}
&W(f \otimes f_1 \otimes f_2 \otimes \ldots \otimes f_n) \\
&\quad= X(f) W(f_1 \otimes f_2 \otimes \ldots \otimes f_n) 
- \sum_{i=1}^n q^{i-1} \Exp{f f_i} W(f_1 \otimes \ldots \otimes \check{f}_i \otimes \ldots \otimes f_n) \\
&\quad\quad- \sum_{i=1}^n q^{i-1} W(f f_i \otimes \ldots \otimes \check{f}_i \otimes \ldots \otimes f_n)
- \Exp{f} W(f_1 \otimes f_2 \otimes \ldots \otimes f_n).
\end{split}
\end{equation}

\begin{Lemma}
For $0 \leq q < 1$,
\[
\norm{X(f)} \leq \Bigl( 1 + \frac{1}{\sqrt{1-q}} \Bigr)^2 \norm{f}.
\]
\end{Lemma}

\begin{proof}
Since both states $\phi$ and $\Exp{\cdot}$ are faithful,
\begin{align*}
\norm{X(f)}
&= \lim_{n \rightarrow \infty} \norm{X(f)}_{L^n(\phi)}
= \Bigl( \sum_{\pi \in \Part(n)} q^{\rc{\pi}} \prod_{B \in \pi} \Exp{f^{\abs{B}}} \Bigr)^{1/n} \\
&\leq \Bigl( \sum_{\pi \in \Part(n)} q^{\rc{\pi}} \prod_{B \in \pi} \norm{f}^{\abs{B}} \Bigr)^{1/n}
= \Bigl( \sum_{\pi \in \Part(n)} q^{\rc{\pi}} \Bigr)^{1/n} \norm{f}
= \norm{X(\id)}_n \norm{f}.
\end{align*}
Thus
\[
\norm{X(f)} \leq \norm{X(\id)} \cdot \norm{f}.
\]
The distribution of $X(\id)$ is the standard $q$-Poisson distribution (the orthogonality measure of the continuous big $q$-Hermite polynomials). The distribution, and in particular its support, are known explicitly, see \cite{SaiKraw} and their references. Its maximum is the indicated constant.
\end{proof}

The following are some analogs of formulas of the main body of the paper for this context, mostly given without proof. The Kailath-Segall formula takes the form
\begin{equation*}
W \left(f^{\otimes (n+1)} \right)
= \sum_{k=0}^n (-1)^k \frac{[n]_q!}{[n-k]_q!} X(f^{k+1}) W \left(f^{\otimes (n-k)} \right)
- \Exp{f} W \left(f^{\otimes n} \right).
\end{equation*}

\begin{Prop}
\begin{equation}
\label{Wick-product}
X(f_1) X(f_2) \ldots X(f_n) =
\sum_{\pi \in \Part(n)} \sum_{S \subset \pi} q^{\rc{S, \pi}} \prod_{B \not \in S} \Exp{\prod_{i \in B} f_i} W \left(\bigotimes_{B \in S} \prod_{i \in B} f_i \right) 
\end{equation}
\end{Prop}

\begin{proof}
Evaluate both sides of the expression~\eqref{Wick-product} on $\Omega$. We obtain
\[
\prod_{i=1}^n \bigl( a(f_i) + a^\ast(f_i) + p(f_i) + \Exp{f_i} \bigr) \Omega = \sum_{\pi \in \Part(n)} \sum_{S \subset \pi} q^{\rc{S, \pi}} \prod_{B \not \in S} \Exp{\prod_{i \in B} f_i} \bigotimes_{B \in S} \biggl( \prod_{i \in B} f_i \biggr).
\]
Each term on the right-hand-side, possibly up to a power of $q$, is obtained by applying a unique sequence of operators $Z_1 Z_1 \ldots Z_n$ to $\Omega$, where
\[
Z_i =
\begin{cases}
a(f_i) & \text{if } i \in B, B \not \in S, i = \min(B), \abs{B} > 1, \\
a^\ast(f_i) & \text{if } i \in B, i = \max(B), \text{ either } B \in S \text{ or } (B \not \in S, \abs{B} > 1), \\
p(f_i) & \text{if } i \in B, i \neq \max(B), \text{ either } B \in S \text{ or } (B \not \in S, i \neq \min(B)), \\
\Exp{f_i} & \text{if } i \in B, B \not \in S, \abs{B} = 1.  
\end{cases}
\]
The rest of the proof proceeds as in Proposition~\ref{Prop:Product-Wick}.
\end{proof}

In particular, since $\state{W(\cdot)} = 0$,
\begin{equation*}
\state{X(f_1) X(f_2) \ldots X(f_n)} =
\sum_{\pi \in \Part(n)} q^{\rc{\pi}} \prod_{B \in \pi} \Exp{\prod_{i \in B} f_i},
\end{equation*}
a generalization of equation~\eqref{Moment}. See \cite{AnsAppell} for further results in this direction.

Now let
\[
H = L^2(\mf{R}_+, dt) \otimes L^2(\mf{R}, \nu),
\]
where $\nu$ is a compactly supported probability measure. Then polynomials are contained in and are dense in $L^2(\mf{R}, \nu)$. Denote 
\[
X(t) = X(\chf{[0, t)} \otimes x)
\]
and
\[
\Delta_k(t) = X(\chf{[0, t)} \otimes x^k).
\]
Note that we no longer have a shift in the index. As in Lemma~\ref{Lemma:Diagonal}, 
\[
\Delta_k(t) = \int_0^t (dX(t))^k = \lim_{\delta(\mc{I}) \rightarrow 0} \sum_{i=1}^N X(\chf{I_i} \otimes x)^k
\]
exists in $L^2(\Gamma(\mc{A}), \varphi)$, where $\mc{I} = \set{I_i}_{i=1}^N$ is a subdivision of $[0, t)$. Consequently, $\Omega$ is cyclic and separating for the algebra generated by$\set{X(t): t \in \mf{R}_+}$.

Define $\psi_n(t)$ and $\St{\pi}(t)$ as in Definition~\ref{Defn:St}. Then 
\begin{equation*}
\St{\pi}(t)
= \sum_{S \subset \pi} q^{\rc{S, \pi}} \left( \prod_{B \not \in S} t \Exp{x^{\abs{B}}} \right)
W \left(\chf{[0, t)^{\abs{S}}} \otimes \bigotimes_{B \in S} x^{\abs{B}} \right).
\end{equation*}
So the Wick product decomposition \eqref{Wick-product} is just the elementary combinatorial decomposition
\[
X(t)^n = \sum_{\pi \in \Part(n)} \St{\pi}(t).
\]

Finally, assume $\Exp{x} = \int_{\mf{R}} x \,d\nu(x) = 0$. Then
\[
\psi_n(t) = W((\chf{[0, t)} \otimes x)^{\otimes n}).
\]

More generally, for $F \in L^2(\mf{R}_+^n, dt^{\otimes n})$ and 
\[
Y_k(t) 
= W(\chf{[0, t)} \otimes x^k) 
= \Delta_k(t) - t \Exp{x^k},
\]
\begin{equation*}
\int F(t_1, \ldots, t_n) dY_{u(1)}(t_1) \ldots dY_{u(n)}(t_n)
= W(F \otimes (x_1^{u(1)} x_2^{u(2)} \ldots x_n^{u(n)})).
\end{equation*}


\providecommand{\bysame}{\leavevmode\hbox to3em{\hrulefill}\thinspace}
\providecommand{\MR}{\relax\ifhmode\unskip\space\fi MR }
\providecommand{\MRhref}[2]{%
  \href{http://www.ams.org/mathscinet-getitem?mr=#1}{#2}
}
\providecommand{\href}[2]{#2}

\end{document}